\documentclass{au}
\usepackage{amssymb,latexsym}
\usepackage{newlattice}

\DeclareMathOperator{\Princ}{Princ}

\DeclareMathOperator{\Lat}{Lat}

\DeclareMathOperator{\Frame}{Frame}
\DeclareMathOperator{\Frucht}{Frucht}

\newtheorem{theorem}{Theorem}

\theoremstyle{definition}

\begin{document}
\title[The Independence Theorem]
{Homomorphisms and principal congruences\\ of bounded lattices. III. 
\\The Independence Theorem} 
\author{G. Gr\"{a}tzer} 
\email[G. Gr\"atzer]{gratzer@me.com}
\address{Department of Mathematics\\
  University of Manitoba\\
  Winnipeg, MB R3T 2N2\\
  Canada}
\urladdr[G. Gr\"atzer]{http://server.maths.umanitoba.ca/homepages/gratzer/}
\date{April 13, 2016}

\begin{abstract}
A new result of G. Cz\'edli states that 
for an ordered set $P$ with at least two elements 
and a group $G$, there exists a bounded lattice $L$ such 
that the ordered set of principal congruences of $L$ 
is isomorphic to $P$ 
and the automorphism group of~$L$ is isomorphic to $G$.

I provide an alternative proof utilizing a result of mine
with J. Sichler from the late 1960-s.
\end{abstract}

\maketitle

\section{Introduction}\label{S:Introduction}
For a bounded lattice $L$, let $\Princ L$ denote
the ordered set of principal congruences of $L$. 
I proved in~\cite{gG14} that $\Princ L$ can be characterized
as a bounded ordered set.\footnote{See the references 
for all the papers that built on this result.}
Let $\Aut L$ denote the group of automorphisms
of a lattice $L$. 
G.~Cz\'edli~\cite{gC16a} proved that these two structures are independent.

\begin{theorem}[Independence Theorem]\label{T:main}
Let $P$ be a bounded ordered set with at least two elements 
and let $G$~be a group. 
Then there exists a bounded lattice $L$ 
such that $\Princ L$ is order isomorphic to $P$
and $\Aut L$ is group isomorphic to~$G$. 
\end{theorem}

I will sketch a new approach to this result.

Note that Cz\'edli proved more than what is stated in Theorem~\ref{T:main}.  
He constructed a lattice $L$ that is selfdual and of length $16$. 
To prove this, he had to work much harder than I do in this note.

\subsection*{Notation}
We use the notation as in \cite{CFL2}.
You can find the complete

\emph{Part I. A Brief Introduction to Lattices} and  
\emph{Glossary of Notation}

\noindent of \cite{CFL2} at 

\verb+tinyurl.com/lattices101+

\section{Families of rigid lattices}
\label{S:rigid}
We call a family $\setm{K_i}{i \in I}$ 
of bounded lattices \emph{mutually rigid} if $K_i$ has no 
$\set{0,1}$-embedding into $K_j$ for $i \neq j \in I$. 
Of course, all $K_i$, $i \in I$ are rigid---they only have the trivial automorphism.

Our proof of the Independence Theorem is based on the following result,
see G.~ Gr\"atzer and J. Sichler \cite{GS70}.

\begin{theorem}\label{T:GrSi}
For every cardinal $\fm$,
there exists a mutually rigid family $\setm{K_i}{i \in I}$ of simple
bounded lattices satisfying $|I| = \fm$.
\end{theorem}

Theorem~\ref{T:GrSi} is a very special case of the result 
of G.~ Gr\"atzer and J. Sichler \cite{GS70}, 
which deals with endomorphism semigroups. 
However, Theorem~\ref{T:GrSi} adds ``simple'', and we will now justify it.

Let $G$ be an infinite connected graph 
in which every element is in some cycle of odd length. 
Let $\Lat G$ denote the free bounded lattice generated by $G$
with respect to the condition that all edges $\set{x,y}$ of $G$
be complemented pairs in $\Lat G$. It follows from  Corollary 567 in LTF, \cite{LTF} 
(originally, C.\,C. Chen and G.~Gr\"atzer \cite{CG69b} and G.~Gr\"atzer \cite{gG71}),
that a pair of elements $\set{x,y} \neq \set{0,1}$
is complementary in $\Lat G$ if{}f $\set{x,y}$ is an edge of~$G$.
Let $G$ be a rigid graph. In G.~ Gr\"atzer and J. Sichler \cite{GS70}, 
we proved that $\Lat G$ is a rigid lattice.

Now we extend $\Lat G$ to a lattice $\Lat^+ G$.
For every $a \in L$ with $0 < a$, we add an atom $p_a < a$
so that 
\begin{align}
    p_a \jj x &= a,\label{E:join}\\
    p_a \mm x &= 0\label{E:meet}
\end{align}
for all $0 < x < a$.
We~observe that the new lattice $\Lat^+ G$ is rigid,
because we can recognize the elements of $G$ 
as complemented elements that are not atoms or dual atoms.

We claim that $\Lat^+ G$ is simple.
Indeed, if $\bga > \zero$ is a congruence of $\Lat^+ G$, then there are 
$u < v \in \Lat^+ G$ such that 
\begin{equation}\label{E:cong}
   \cng u=v(\bga).
\end{equation}
By \eqref{E:meet},  
\begin{equation}\label{E:cong2}
   \cng 0=p_v(\bga).
\end{equation}
It follows by \eqref{E:join} that 
\begin{equation}\label{E:cong3}
   \cng 1=p_1(\bga)
\end{equation}
and therefore,
\begin{equation}\label{E:cong4}
   \cng 0=x(\bga)
\end{equation}
for all $x \in \Lat^+ G - \set{1, p_1}$.

Take an edge $\set{x_1, x_2}$ of $G$. Applying \eqref{E:cong4} twice,
we obtain that 
\begin{equation}\label{E:cong5}
   \cng 0=x_i(\bga)
\end{equation}
for $i = 1,2$. Joining these two congruences, we conclude that
\begin{equation}\label{E:cong6}
   \cng 0=1(\bga),
\end{equation}
that is $\bga = \one$, and so $\Lat^+ G$ is simple.

So $\Lat^+ G$ is simple and rigid.

By Z.~Hedrl\'\i n and A. Pultr \cite{HP64}, 
we can take a family $\setm{G_i}{i \in I}$ 
of mutually rigid, connected graphs with $|I| = \fm$, 
in which every element is in some cycle of odd length.
Then the family $\setm{\Lat^+G_i}{i \in I}$ is rigid and simple; in fact,
it is mutually rigid for the same reason as we argued before.

\section{Proving the Independence Theorem}\label{S:Independence}
Let $P$ be the order in Theorem~\ref{T:main}. 
Let $0$ and $1$ denote the zero and unit of~$P$, respectively.

We start out by constructing the lattice $\Frame P$, as in \cite{gG14},
see also \cite{CFL2}.

\begin{figure}[h!t]
\centerline{\includegraphics[scale=1.0]{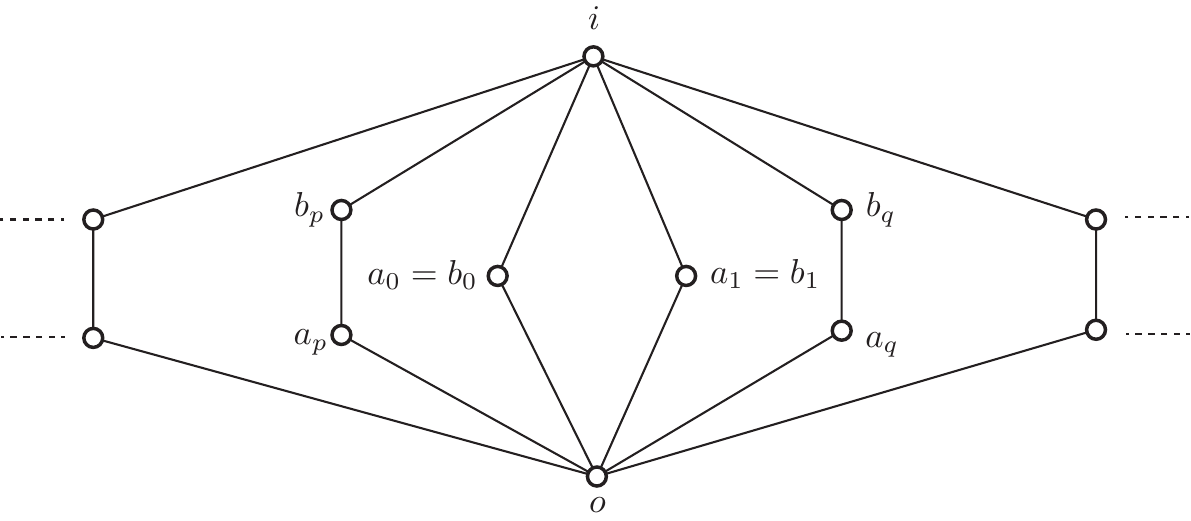}}
\caption{The lattice $\Frame P$}\label{F:Frame}
\end{figure}

Let $\Frame P$ consist of the elements \text{$o$, $i$} 
and the elements $a_p, b_p$ for every $p \in P$,
where $a_p \neq b_p$ for every $p \in P^-$
and $a_0 = b_0$, $a_1 = b_1$. These elements are ordered and
the lattice operations are formed as in Figure~\ref{F:Frame}. 

We are going to construct a lattice $K$ as an extension of $\Frame P$. 
\begin{figure}[hbt]
\centerline{\includegraphics[scale=1.0]{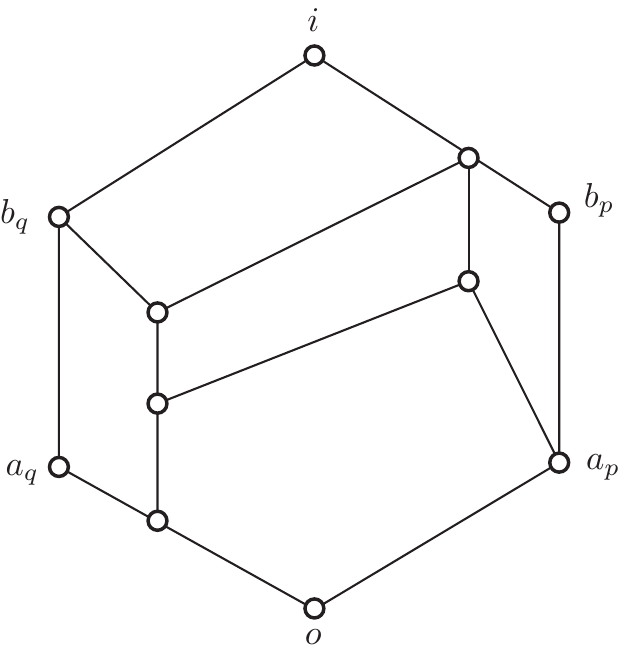}}
\caption{The sublattice}\label{F:S} 
\end{figure}
For $p < q \in P^-$, we add five elements to the sublattice 
$
   \set{o, a_p, b_p, a_q, b_q, i}
$
of $\Frame P$ as illustrated in Figure \ref{F:S}.
In $K$, the principal congruence $\con{a_p, b_p}$ 
represents $p \in P^-$.

Now we prove Theorem~\ref{T:main}. 

Let $P$ and $G$ be given as in this theorem. 
The lattice $K$ we have just constructed satisfies that $\Princ K \iso P$.
However, $\Aut K \iso \Aut P$.

Let $I = P^-$. For every $[o, a_p]$, we insert the lattice $K_p$,
provided by Theorem~\ref{T:GrSi},  
into $[o, a_p] \ci K$ identifying the bounds,
obtaining the lattice $\ol K$. 
We then have $\Princ \ol K \iso P$ and $\Aut \ol K$ is rigid.

Now take the Frucht lattice, $\Frucht G$
(see \cite{CFL2}, \cite{LTF}, R.~Frucht~\cite{rF38} and \cite{rF50}); 
it is a lattice of length $3$, satisfying $\Aut(\Frucht G) \iso G$.

\begin{figure}[h!t]
\centerline{\includegraphics[scale=0.7]{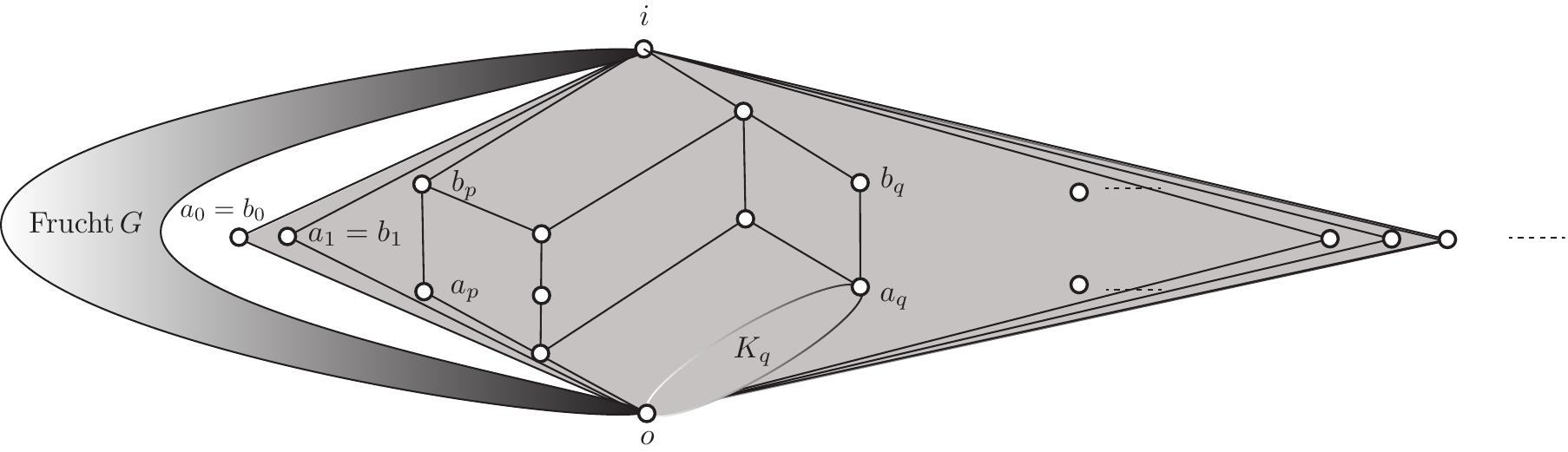}}
\caption{The final step}\label{F:Framenew}
\end{figure}

We obtain the lattice $L$ by forming the disjoint union of $\ol K$
and $\Frucht G$ and identifying the bounds, see Figure~\ref{F:Framenew}.
The diagram does not show $K_p$.

\newpage

\end{document}